\documentclass{IEEEtran4PSCC}
\ifCLASSINFOpdf
   \usepackage[pdftex]{graphicx}
\else
   \usepackage[dvips]{graphicx}
\fi
\usepackage[cmex10]{amsmath}
\usepackage{amsthm}
\usepackage{nccmath}
\usepackage{comment}
\usepackage{eurosym}
\usepackage{xcolor}

\usepackage{url}

\usepackage{nomencl}
\usepackage{ifthen}
\usepackage{amssymb}
\usepackage{threeparttable}
\usepackage{booktabs} 
\usepackage{siunitx} 
\usepackage{array} 

\makeatletter
\let\old@ps@headings\ps@headings
\let\old@ps@IEEEtitlepagestyle\ps@IEEEtitlepagestyle
\def\psccfooter#1{%
    \def\ps@headings{%
        \old@ps@headings%
        \def\@oddfoot{\strut\hfill#1\hfill\strut}%
        \def\@evenfoot{\strut\hfill#1\hfill\strut}%
    }%
    \def\ps@IEEEtitlepagestyle{%
        \old@ps@IEEEtitlepagestyle%
        \def\@oddfoot{\strut\hfill#1\hfill\strut}%
        \def\@evenfoot{\strut\hfill#1\hfill\strut}%
    }%
    \ps@headings%
}
\makeatother

\psccfooter{%
        \parbox{\textwidth}{\hrulefill \\ \small{23rd Power Systems Computation Conference} \hfill \begin{minipage}{0.2\textwidth}\centering \vspace*{4pt} \includegraphics[scale=0.06]{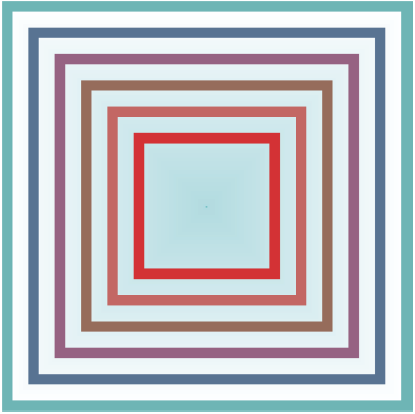}\\\small{PSCC 2024} \end{minipage} \hfill \small{Paris, France --- June 4 -- 7, 2024}}%
}

\usepackage{etoolbox}
\renewcommand{\nomgroup}[1]{%
\ifthenelse{\equal{#1}{S}}{\item[\textbf{Sets}]}{%
\ifthenelse{\equal{#1}{P}}{\item[\textbf{Parameters}]}{%
\ifthenelse{\equal{#1}{V}}{\item[\textbf{Variables}]}{}}}
}
\makenomenclature

\begin{document}
%
\title{An Incentive Regulation Approach for Balancing Stakeholder Interests in Transmission Merchant Investment}
\color{black}

\author{
\IEEEauthorblockN{Yuxin Xia}
\IEEEauthorblockA{School of Engineering \\
The University of Edinburgh\\
Edinburgh, U.K.\\
yuxin.xia@ed.ac.uk}
\and
\IEEEauthorblockN{Iacopo Savelli}
\IEEEauthorblockA{The Centre for Research on Geography, \\Resources, Environment, Energy \& Networks\\
Bocconi University\\
20100 Milano, MI, Italy\\
iacopo.savelli@unibocconi.it}
\and
\IEEEauthorblockN{Thomas Morstyn}
\IEEEauthorblockA{Department of Engineering Science\\ University of Oxford\\
Oxford, U.K.\\
thomas.morstyn@eng.ox.ac.uk}
}

\maketitle

\begin{abstract}
The merchant-regulatory mechanism represents a promising tool that combines the benefits of merchant investment and regulated investment, thereby providing efficient incentives for merchant Transmission Companies (Transcos) subject to regulatory compliance. Taking the H-R-G-V mechanism as a foundational example of this approach, it permits Transcos to receive the total surplus increase from investments, and the profit-maximizing Transco will perform social welfare maximum investment under this mechanism. However, one drawback of this mechanism is that it allows the Transco to receive the whole benefit created by the Transco, while excluding consumers and generators from the resultant economic benefits. To address this issue, we propose an incentive tuning parameter, which is incorporated into the calculation of the incentive fee for the Transco. Accordingly, the regulatory framework can effectively manage the Transco's profit and allow market participants to access economic benefits, thus ensuring a fair distribution of economic advantages among the stakeholders, while the impact on overall social welfare remains relatively modest. The results on the case study demonstrate that this careful balancing act maintains the essence of the H-R-G-V mechanism while addressing its critical gap---the equitable sharing of economic gains.

\end{abstract}
\color{black}

\thanksto{\noindent Submitted to the 23rd Power Systems Computation Conference (PSCC 2024).}

\section{Introduction}
Transmission investment is essential to the success of energy transition, providing consumers with non-discriminatory access to affordable generation and ensuring competitiveness and sustainability \cite{Teusch2012}. According to the report in \cite{ukreport}, the investment in transmission and distribution grids is expected to increase to \euro 40-62 billion per year in the EU to meet climate and energy goals. In the US, the total capital investment in transmission network is expected to reach 3.7 trillion dollars by 2050 \cite{usreport}. Therefore, it is important to deliver efficient investments and relieve congestion problems, while ensuring fairness between stakeholders \cite{WilliamHogan}.

Historically, electric transmission has been regarded as a natural monopoly due to its inherent characteristics \cite{ROSELLON20113}. In contrast to this centralized approach, technological advancements in the transmission sector have given rise to a decentralized, market-driven model known as merchant transmission investment \cite{BOFFA2015455,b68c23ba-54e0-37c9-b945-aa473f4789f0}. \color{black} Unlike regulated monopoly transmission investment, merchant investment fosters a market-driven environment and encourages unrestricted competition in the investment process \cite{JoskowPaul2002}. Merchant investors seek remuneration through the the sale of financial or physical transmission rights, or the congestion rent \cite{Joskow2020,KRISTIANSEN20104107,Hogan2002FINANCIALTR}. Nevertheless, despite the numerous advantages associated with merchant transmission investment (as discussed in \cite{leautier2023}), significant concerns arise regarding its potential to lead to sub-optimal expansion and questions remain related to both theoretical design and real-world implementation \cite{8322275}.

Despite the extensive reform in the electricity industry, the transmission sector of electric power systems has largely remained under regulation \cite{Majidi2020}. In many countries, such as in England and Wales \cite{JoskowPaul2002}, transmission companies continue to function as natural monopolies, necessitating regulatory incentives to promote investment \cite{Majidi2020}. The regulatory framework entails various design approaches, including cost-of-service mechanisms \cite{Majidi2020}, price-cap regulation mechanisms \cite{Vogelsang2001,cite-key},  incentive regulation approaches \cite{ARMSTRONG20071557} and merchant-regulated mechanisms \cite{Vogelsang2020,Khastieva2021,Hogan:2010uf,Hesamzadeh2018}.

According to the argument in \cite{Joskow2005}, merchant investment is considered to be a supplementary approach rather than a substitute for regulatory investment. In fact, Australia utilizes this combination of merchant and regulated investment strategies \cite{JoskowPaul2002}. This merchant-regulated mechanism combines the benefits of both merchant and regulated investment, with the H-R-G-V mechanism serving as an example \cite{Vogelsang2020,Khastieva2021,Hesamzadeh2018,Disjunctive}. This mechanism builds upon the price-cap theory \cite{Vogelsang2001} and the Incremental Surplus Subsidy (ISS) scheme \cite{Sappington1988}. By determining a regulated incentive fee that depends on the total contribution to economic benefits from their investment, the H-R-G-V incentive mechanism aims to provide incentives to the profit-maximizing Transco for performing social-welfare maximizing investments\footnote{Readers are referred to Section \ref{section2a} for a comprehensive understanding of the H-R-G-V mechanism.}.

\color{black}
The primary objective of the Transco typically revolves around profit maximization, while the regulator is entrusted with the responsibility of promoting social welfare \cite{Vogelsang2020}. Moreover, a ‘benevolent' regulator can also establish an objective function that assigns weights to either consumer benefits or the net profit of the Transco (the regulated firm) \cite{ARMSTRONG20071557}. Concerns have been raised about the H-R-G-V mechanism, which allows Transcos to capture \emph{all} the welfare improvements. This puts consumers and generators at a disadvantage, as they do not receive any economic benefits from the network expansion as the incentive fee is extracted from them. This situation raises a significant research question: \emph{How to develop a regulated mechanism for the Transco that not only effectively addresses the interests of both public and private entities but also takes into account the benefits for market participants?} To address this issue, we propose an incentive scheme tuning parameter. By limiting the income of the Transco through this parameter, regulators can explore the trade-offs between the Transco's profits, social welfare, and the interests of market participants. Through appropriate adjustments of this parameter, regulators can strike a balance between incentivizing efficient investment, achieving modest social welfare and ensuring an equitable distribution of economic benefits among all stakeholders. 

The rest of the paper is organized as follows: Section~\ref{section2} introduces the merchant-regulatory mechanism and presents the bi-level optimization problem. Section~\ref{section3} presents the reformulation of the bilevel problem as a mixed-integer linear
program (MILP) problem. Results for two case studies are discussed in Section~\ref{section4}. Lastly, Section~\ref{section5} concludes the paper.

\section{Problem Formulation}\label{section2}
This section specifies the Transco's investment problem under the proposed merchant-regulatory incentive mechanism. Section~\ref{section2a} introduces the incentive tuning parameter and the associated incentive mechanism. The bilevel model is described in Section~\ref{section2b} and the math formulation of the bilevel optimization problem is specified in Section~\ref{section2c} and~\ref{section2d}.

\subsection{The merchant-regulatory mechanism with the incentive tuning parameter}\label{section2a}
\begin{figure}[tp]
    \centering    \includegraphics[width=\columnwidth]{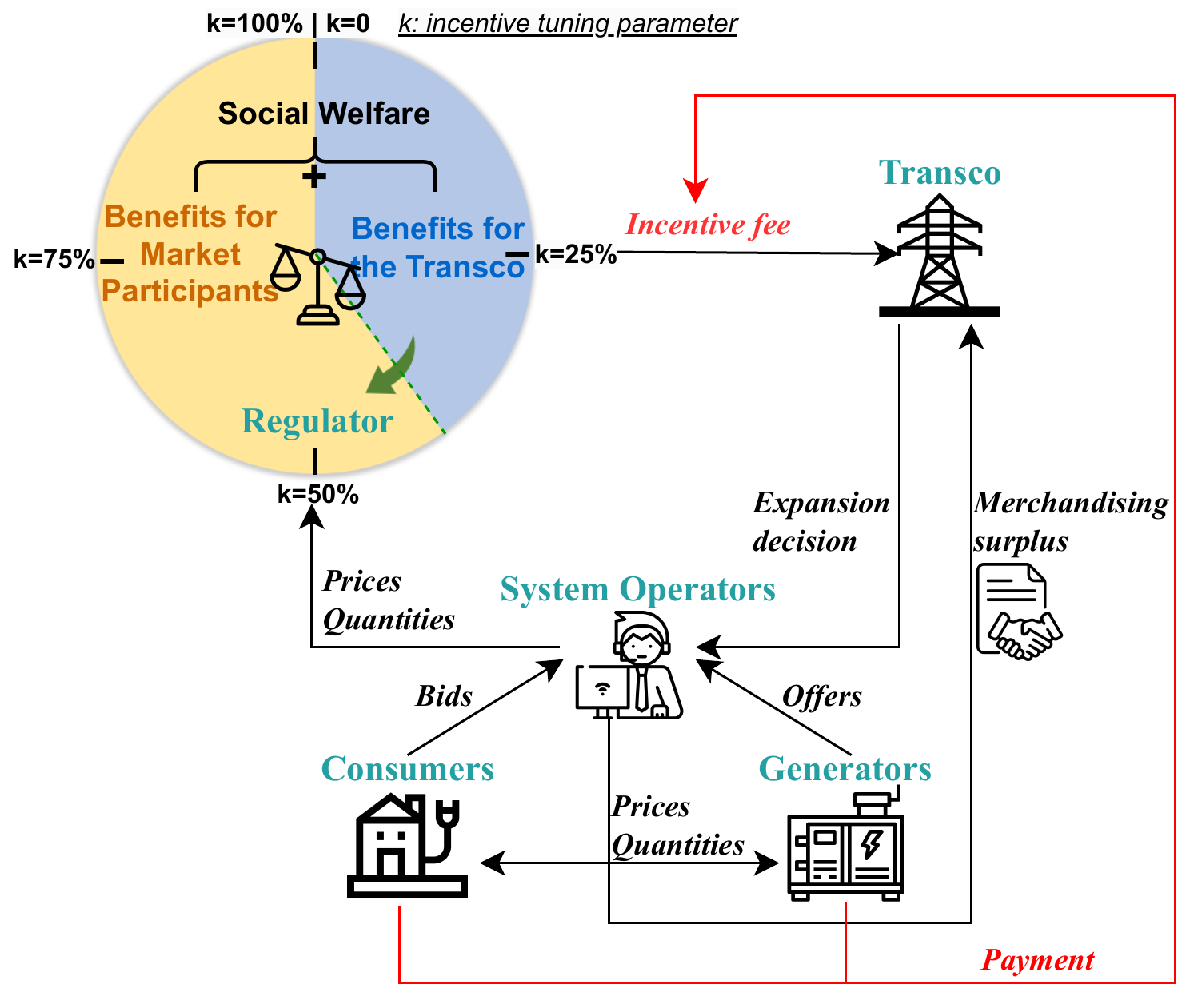}
    \caption{The proposed regulatory framework for managing transmission investments, highlighting the interaction between the Transcos, regulators, SOs, and market participants. The figure emphasizes the role of the incentive tuning parameter, $\kappa$, set by regulators, which determines the allocation of economic benefits between market participants and the Transco. A higher $\kappa$ value indicates a larger proportion of benefits directed towards the Transco. The figure also shows the concept of social welfare, representing the cumulative benefits to both the Transco and market participants.}
    \label{fig:regulation}
\end{figure}
Fig.~\ref{fig:regulation} illustrates the regulatory framework governing the interactions among the regulated Transco, market participants, regulators, and System Operators (SOs). Collaboration and information sharing form the basis of the regulated incentive mechanism within this framework. The Transco is responsible for planning and bears the costs of investment. SOs not only operate and optimize the market based on participants' bids to maximize social welfare but also facilitate electricity transactions and compute the merchandising surplus for the Transco. The regulator, typically a government body, receives dispatch information from the SOs and calculates the regulatory incentive fee to be paid by generators and consumers in an aggregated form, which is subsequently paid to the Transco. When determining the incentive fee, the regulator considers three key factors through incentive tuning parameter adjustments: social welfare (representing public interests, which is also the sum of the benefits for the Transco and market partcipants), the Transco's profit (representing private interests), and the benefits received by producers and consumers (representing market participants' interests). 

We denote the line expansion decisions, line expansion cost, merchandising surplus, incentive fee, generator surplus and load surplus at investment planning period $t$ as $u_t$, $C_t(u_t)$, $MS_t$, $\Phi_t$, $S^{G}_{t}$ and $S^{L}_{t}$, respectively. The discount rate is denoted as $r$.
Social welfare (SW) is calculated as the sum of load surplus $S^{L}_{t}$, generator surplus $S^{G}_{t}$, and merchandising surplus $MS_t$, with the investment cost in new lines $C_t(u_t)$ subtracted from it over the planning horizon, i.e., $SW = \sum_{t \in \mathcal{T}} \frac{1}{(1+r)^{t-1}} \big(MS_t+S^{L}_{t}+S^{G}_{t}-C_t(u_t)\big)$. Furthermore, the Transco's profit $TP_t$ is determined by the sum of merchandising surplus $MS_t$ and the incentive fee $\Phi_t$, with the deduction of line investment costs $C_t(u_t)$, i.e., $TP = \sum_{t \in \mathcal{T}} \frac{1}{(1+r)^{t-1}}\big( MS_t+\Phi_t-C_t(u_t)\big)$. 

In this paper, the incentive fee $\Phi_t$ is calculated based on the incentive tuning parameter $\kappa$ multiplied by the increase in surplus resulting from the investment, where $\kappa \in[0,1]$. The remaining economic benefits resulting from the transmission network investments belong to consumers and generators. In essence, social welfare is the collective benefits received by both the Transco and market participants. Under the modified H-R-G-V mechanism with the incentive tuning parameter $\kappa$, the incentive fee in year $t$ is calculated as 
\begin{align}
\Phi_t=\Phi_{t-1} + \kappa(\Delta S^{G}_{t} +\Delta S^{L}_{t})\label{eq:1}
\end{align}
where $\Delta S^{G}_{t} = S^{G}_{t}-S^{G}_{t-1}$ is the change in the generation surplus, $\Delta S^{L}_{t} =  S^{L}_{t} - S^{L}_{t-1}$ is the change in the load surplus from year $t-1$ to $t$. It is assumed that no investment is performed at $t=1$ and $\Phi_{t=1}=0$. In \cite{Disjunctive}, it is demonstrated that under the original H-R-G-V mechanism, wherein $\kappa = 100\%$, the profit-maximizing investment strategy adopted by Transco aligns with the goal of social welfare maximization. This indicates that Transco, upon receiving the entirety of the surplus increase (economic benefits) resulting from its investments, will engage in investments that optimize social welfare.

The incentive fee, $\Phi_t$, is determined by the regulator using equation~\eqref{eq:1}, without directly specifying the regulator's objective function. As a result, the regulatory constraint expressed in equation~\eqref{eq:1} can be integrated into Transco's upper-level problem for profit calculation, as discussed in \cite{Khastieva2021}. Subsequently, the term ‘\emph{Change in Surplus due to investment}’ is defined as the net increase in generator and load surplus from the first year to any subsequent year $t$, for all $t \in \{\mathcal{T}\backslash 1\}$, represented mathematically as $\sum_{t\in \{\mathcal{T}\backslash 1\}}(S^{G}_{t}+S^{L}_{t}-S^{G}_{1}-S^{L}_{1})$. Based on equation~\eqref{eq:1}, the total incentive fee can be calculated as $\sum_{t\in\mathcal{T}}\Phi_{t} =  \sum_{t\in \{\mathcal{T}\backslash 1\}} \kappa(S^{G}_{t}+S^{L}_{t}-S^{G}_{1}-S^{L}_{1})$, indicating that a $\kappa$ proportion of the ‘Change in Surplus due to investment’ is allocated to Transco, with the remainder distributed among market participants.

\color{black}
The tuning of the incentive parameter $\kappa$ plays a crucial role in managing the Transco's profitability. Decreasing the tuning parameter can enhance the benefits for market participants. However, decreasing this parameter may have adverse consequences for the Transco's investment incentives in the network, consequently compromising overall social welfare. Therefore, the total surplus increase resulting from investment and the benefits for market participants may decrease due to the lack of investment incentives. Consequently, regulators face the challenge of skillfully balancing social welfare, the financial outcomes of the Transco, and the benefits for market participants to achieve an ideal outcome.

The following assumptions are made in this paper:
\begin{itemize}
  
\item We explore the dynamics of network expansion investments, particularly underlining their tendency to be `\emph{lumpy}' as highlighted in previous research \cite{b68c23ba-54e0-37c9-b945-aa473f4789f0}. The notion of `\emph{lumpy expansion}' suggests that capacity expansion in power lines are constrained to predefined, discrete increments rather than continuous scales \cite{SAVELLI2021102450,SAVELLI2020113979}.

\item We assume that both the Transco and the regulatory authority receive all the information of market outcomes \cite{Hesamzadeh2018,Disjunctive}. The incentive fee, denoted as $\Phi_t$, is calculated based on the actual realized surplus, with the regulator wielding the authority to modulate the incentive tuning parameter $\kappa$. It is assumed that the Transco agree on the incentive fee determined by the regulator \cite{Khastieva2021}.

\item The Transco is assumed to be a regulated risk-neutral entity, aiming to maximize profits while bearing the costs associated with line expansion. Its revenue streams are derived from the merchandising surplus and the aforementioned incentive fee \cite{Disjunctive}.

\item We assume a perfect competition among market participants, accommodating the evolving dynamics where demand elasticity is recognized. We assume growing consumer flexibility and price-sensitive behavior, thereby incorporating bids and offers from price-elastic agents, including generators and consumers \cite{SAVELLI2020113979}.

\item Furthermore, we assume that the decision to expand line capacity is made only once within the planning horizon and is characterized by its irreversibility.

\end{itemize}
\color{black}
\subsection{The bilevel optimization model}\label{section2b}
\begin{figure}[thp]
    \centering    \includegraphics[width=\columnwidth]{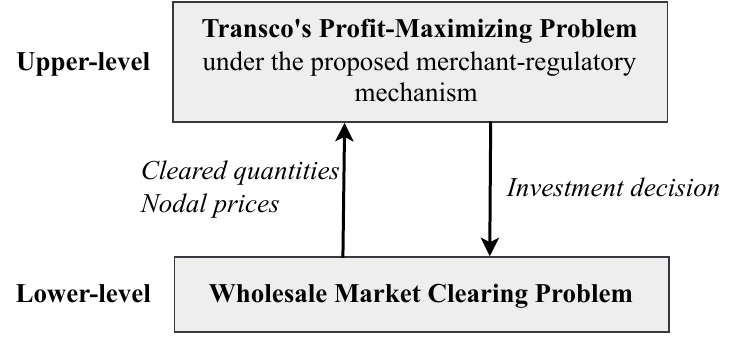}
    \caption{High-level block diagram for the proposed bi-level framework.}
    \label{fig:bilevel}
\end{figure}

In the proposed bilevel optimization model, as depicted in Fig.~\ref{fig:bilevel}, the upper-level problem represents the objective of the profit-maximizing Transco, where the incentive fee is subject to the regulatory constraint~\eqref{eq:1}. The lower-level problem is a standard wholesale market (WSM) clearing problem. The upper-level problem determines the incentive fee $\Phi_t$, binary decision variables $u_{t,l}$ and selected lumpy expansion decision $b^F_{t,l,j}$. Having fixed these upper-level variables, the lower-level problem performs the WSM clearing and determines the allocated quantity for generators $g_{t,s,k,b}$ and consumers $d_{t,s,k,b}$, and WSM prices $\pi_{t,s,b}$.

\subsection{Upper-level problem: Transco's profit maximizing problem}\label{section2c}
Problem~\eqref{eq:upper} specifies the Transco's profit-maximizing problem under the proposed incentive mechanism with the incentive tuning parameter $\kappa$.

\begin{fleqn}
\begin{subequations}
\begin{align}
    &\max_{u_{t,l},b_{t,l, j}^{F},\Phi_t} \sum_{t \in \mathcal{T}} \frac{1}{(1+r)^{t-1}}\bigg(\sum_{s \in \mathcal{S}} \sum_{b \in \mathcal{B}}\sum_{k \in \Omega_{t, s,b}^{D}}\Psi\pi_{t,s,b}^{*}d_{t, s,k,b}^{*} \notag\\
    &- \sum_{s \in \mathcal{S}} \sum_{ b\in \mathcal{B}}\sum_{k \in \Omega_{t, s,b}^{G}}\Psi\pi_{t,s,b}^{*}g_{t, s,k,b}^{*}+\Phi_t-\Psi\sum_{l \in \mathcal{L}}(u_{t,l} K_{l}^{f i x}\notag\\
    &+K_{l}^{v a r} \sum_{j \in \mathcal{J}} b_{t,l, j}^{F} \bar{F}_{l, j})\bigg)\label{eq:upper_obj} \end{align}
    Subject to:
 \begin{align}
    S_{t}^{L}= \sum_{s \in \mathcal{S}} \sum_{b \in \mathcal{B}}\sum_{k \in \Omega_{t, s,b}^{D}} \big( c_{t,s, k,b}^{d} - \pi_{t,s,b}^{*}\big) d_{t, s,k,b}^{*}, \forall t\in\mathcal{T}  \label{eq:upperb} \end{align}
    \vspace{-10pt}
    \begin{align}
    S_{t}^{G}=  \sum_{s \in \mathcal{S}} \sum_{b \in \mathcal{B}}\sum_{k \in \Omega_{t,s, b}^{G}} \big( \pi_{t,s,b}^{*} -  c_{t,s,k,b}^{g}\big) g_{t,s, k,b}^{*},  \forall t\in\mathcal{T} \label{eq:upperc} \end{align}
    \vspace{-10pt}
    \begin{align}
         \Phi_t-\Phi_{t-1}=\kappa\Psi\Big(S_{t}^{L}-S_{t-1}^{L}+S_{t}^{G}-S_{t-1}^{G}\Big), \quad   t\geq 2 \label{eq:upperd} \end{align}
     \vspace{-10pt}
    \begin{align}   
         \sum_{t \in \mathcal{T}}\sum_{j \in \mathcal{J}}b_{t,l, j}^{F} \leq 1 , \quad   \forall l \in \mathcal{L}  \label{eq:uppere} \end{align}
     \vspace{-10pt}
    \begin{align} 
         u_{t,l}=\sum_{j \in \mathcal{J}}b_{t,l, j}^{F}, \quad    \forall l \in \mathcal{L},  \forall t\in\mathcal{T} 
         \label{eq:upperf} \end{align}
    \vspace{-10pt}
    \begin{align}   
         u_{t,l}\in\{0,1\}, b_{t,l, j}^{F}\in\{0,1\}, \Phi_{t=1}=0, u_{t=1,l}=0 \label{eq:upperg} \end{align}
    \label{eq:upper}
\end{subequations}
\end{fleqn}
The objective function~\eqref{eq:upper_obj} represents the aim of the profit-maximizing Transco, whose profit is calculated based on the revenues from the merchandising surplus and the incentive fee $\Phi_{t}$ and the costs of line expansion. Specifically, the merchandising surplus is calculated from optimal solutions of the lower-level problem which depends on the upper level variables (see equations~\eqref{eq:lowerf}-\eqref{eq:lowerg}), including the cleared quantities for demand $d_{t, s,k,b}^{*}$ and generators $g_{t, s,k,b}^{*}$, and WSM prices $\pi_{t,s,b}^{*}$. The line expansion costs consists of a fixed part $K_{l}^{f i x}$ and a variable part $K_{l}^{v a r}$. The term $\bar{F}_{l, j}$ is the lumpy expansion on line $l$, where $l \in\mathcal{L}$, and the line capacity can be increased only by a finite set of discretized quantities $\bar{\mathcal{F}_l}=\bigcup_{j \in J}\bar{F}_{l, j}$ \cite{SAVELLI2020113979}. The discount rate $r$ is used for calculating the present value. Here two time indices $t$ and $s$ are considered, representing the yearly investment period and the hourly operational period, respectively. The parameter $\Psi$ serves as a scaling factor designed to align the costs and benefits from two different time scales. It represents the number of operational periods per investment period, For example, if we assume a daily operational period $\mathcal{S} = \{1,2,\dots,24\}$ and yearly investment periods, the $\Psi$ is set to 365. \color{black}Constraints~\eqref{eq:upperb} and~\eqref{eq:upperc} calculate the load surplus and generator surplus based on the solutions from the lower-level market clearing problem and used to compute the regulated incentive fee, as shown in equation~\eqref{eq:upperd}. The incentive tuning parameter $\kappa$ is applied to the generator and consumer surplus increase due to investment, as discussed in Section~\ref{section2a}. Equations~\eqref{eq:uppere} and~\eqref{eq:upperf} enforce the line expansion decision is taken place once for all investment periods $t\in\mathcal{T}$ and this decision is irreversible. Equation~\eqref{eq:upperg} ensures in the first year, no expansion is performed and the incentive fee is zero.

\subsection{Lower-level problem: Wholesale market clearing problem}\label{section2d}
The WSM clearing problem is specified in Problem~\eqref{eq:lower}.

\begin{fleqn}
\begin{subequations}
\begin{align}
&( d^{*}_{t,s,k,b} ,g^{*}_{t,s,k,b},f^{*}_{t,s,l},\theta_{t,s,b}^{*},[\pi _{t,s,b}^{*}])=\notag\\
&\arg\max\sum_{t \in \mathcal{T}} \sum_{s \in \mathcal{S}} \sum_{b \in \mathcal{B}}\bigg(\sum_{k \in \Omega_{t, s,b}^{D}} c_{t,s, k,b}^{d} d_{t,s, k,b}-\sum_{k \in \Omega_{t, s,b}^{G}} c_{t,s,k,b}^{g}\notag\\
& g_{t,s, k,b}\bigg)
\label{eq:lower_obj}
\end{align}
Subject to:
\begin{align}
&-\sum_{k\in\Omega_{t, s,b}^{G}}g_{t,s, k,b}+\sum_{k \in \Omega_{t, s,b}^{D}}d_{t, s,k,b}+\sum_{l\in\mathcal{L}}S_{l,b}f_{t,s,l}- \sum_{l\in\mathcal{L}}R_{l,b}\notag\\
&f_{t,s,l}=0, \quad \forall t\in \mathcal{T},\forall s\in \mathcal{S}, \forall b\in \mathcal{B}\hspace{25pt} [\pi_{t,s,b}\in \mathbb{R}] \label{eq:lowerb} \end{align}
\vspace{-15pt}
\begin{align}   
&g_{t, s,k,b}^{ \min }  \leq g_{t, s,k,b} \leq g_{t, s,k,b}^{ \max },\quad\forall t \in \mathcal{T}, \forall s\in \mathcal{S}, \forall k \in \Omega_{t, s,b}^{G},\notag\\ 
&\forall b \in \mathcal{B}\hspace{90pt} [\varphi^{G,\min}_{t,s,k,b}\geq 0,\varphi^{G,\max}_{t,s,k,b}\geq 0]  \label{eq:lowerc} \end{align}
\vspace{-15pt} 
\begin{align}   
&d_{t, s,k,b}^{ \min } \leq d_{t, s,k,b}\leq d_{t, s,k,b}^{ \max },\quad\forall t \in \mathcal{T}, \forall s\in \mathcal{S}, \forall k \in \Omega_{t,s, b}^{D},\notag\\ 
&\forall b \in \mathcal{B}\hspace{87pt} [\varphi^{D,\min}_{t,s,k,b}\geq 0,\varphi^{D,\max}_{t,s,k,b}\geq 0]\label{eq:lowerd} \end{align}
\vspace{-15pt}
\begin{align}   
&f_{t,s,l}=B_{l}\left(\sum_{b\in\mathcal{B}}S_{l,b}\theta_{t,s,b}-\sum_{b\in\mathcal{B}}R_{l,b}\theta_{t,s,b}\right), \notag\\
& \hspace{75pt} \forall t \in \mathcal{T}, \forall s\in \mathcal{S},\forall l \in \mathcal{L} \quad [\gamma_{t,s,l}\in \mathbb{R}] \label{eq:lowere} \end{align}
\vspace{-15pt}
\begin{align} 
&f_{t,s, l}\leq \mathcal{F}^{0}_{l}+\sum_{\hat{t} \in\{2, . . t\}}\sum_{j\in\mathcal{J}}b^F_{\hat{t},l,j}\overline{F}_{l,j},\quad\forall t \in \mathcal{T},\forall s\in \mathcal{S},\forall l \in \mathcal{L}\notag \\
& \hspace{180pt} [\mu^{\max}_{t,s,l}\geq 0]
\label{eq:lowerf} \end{align}
\vspace{-15pt}
\begin{align}   
&-f_{t,s, l}\leq \mathcal{F}^{0}_{l}+\sum_{\hat{t} \in\{2, . . t\}}\sum_{j\in\mathcal{J}}b^F_{\hat{t},l,j}\overline{F}_{l,j}, \forall t \in \mathcal{T},\forall s\in \mathcal{S},\forall l \in \mathcal{L}\notag\\
&\hspace{180pt} [\mu^{\min}_{t,s,l}\geq 0]  \label{eq:lowerg} \end{align}
\vspace{-15pt}
\begin{align}   
&-\theta^{max}_b\leq \theta_{t,s,b} \leq \theta^{max}_b,\quad \forall t \in \mathcal{T},\forall s\in \mathcal{S},\forall b \in \mathcal{B}\notag\\
&\hspace{135pt} [\xi^{\min}_{t,s,b}\geq 0,\xi^{\max}_{t,s,b}\geq 0]  \label{eq:lowerh} \end{align}
\vspace{-15pt}
\begin{align}   
\theta_{t,s,1}=0,\quad \forall t \in \mathcal{T},\forall s\in \mathcal{S} \hspace{65pt} [\chi_{t,s}\in\mathbb{R}] \label{eq:loweri} \end{align}
    \label{eq:lower}
\end{subequations}
\end{fleqn}
The objective of the WSM clearing problem is to maximize social welfare spanning all planning years $t\in\mathcal{T}$ and operation periods $s\in\mathcal{S}$. The power balance equation is modelled in~\eqref{eq:lowerb} in node $b$ at year $t$ and operational period $s$. The supply and demand upper and lower limits of generator and consumers are shown in equations~\eqref{eq:lowerc} and ~\eqref{eq:lowerd}, respectively. The proposed model employs the DC optimal power flow (DC-OPF) approximation and the resulting power flow is modelled in~\eqref{eq:lowere}. The flow limit is enforced by~\eqref{eq:lowerf} and~\eqref{eq:lowerg}. Binary variables $b^F_{\hat{t},l,j}$ are the lumpy expansion decision determined in the upper-level problem. The product $\sum_{\hat{t} \in\{2,...,t\}}\sum_{j\in\mathcal{J}}b^F_{\hat{t},l,j}\overline{F}_{l,j}$ determines the selected amount of lumpy expansion from $\hat{t} \in\{2,...,t\}$ where $t\in\mathcal{T}$ since the investment decision is irreversible. Equations~\eqref{eq:lowerh} and \eqref{eq:loweri} define the range of voltage phase angle of node $b$. 

\section{Reformulation as a MILP problem}\label{section3}
\begin{figure}[tp]
    \centering    \includegraphics[width=\columnwidth]{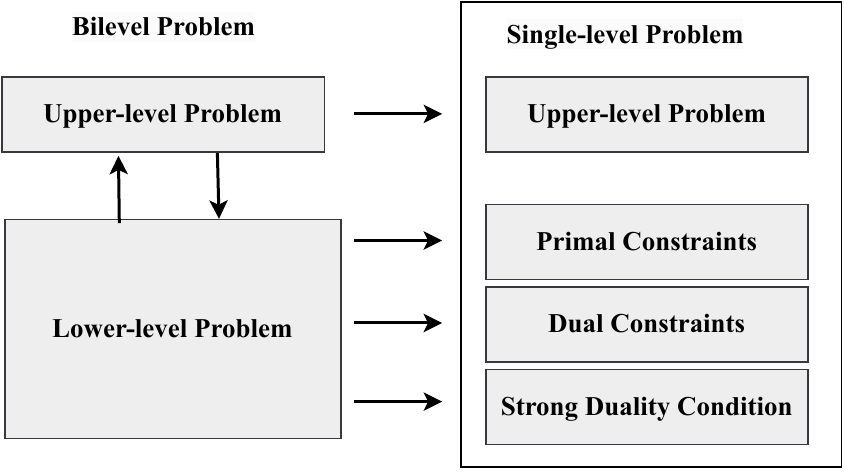}
    \caption{Reformulation scheme.}
    \label{fig:reformulation}
\end{figure}

In this section, the bilevel model (Problem \eqref{eq:upper} and Problem \eqref{eq:lower}) is reformulated into a single-level problem by resorting to the primal constraints (equations~\eqref{eq:lowerb}-\eqref{eq:loweri}), dual constraints (equations~\eqref{eq:WSMdual}) and the strong duality condition (equation~\eqref{eq:strong_duality}) of the lower-level problem \cite{SAVELLI2021102450}, as shown in Fig.~\ref{fig:reformulation}. In addition, the nonlinear terms are removed (discussed in Section~\ref{section3c}) and the final MILP problem is presented in Section~\ref{section3d}.

The bilevel problem stated in Section~\ref{section2c} and Section~\ref{section2d} can be equivalently recast as a single-level problem, which is defined as follows:

\begin{fleqn}
\begin{subequations}
\begin{align}
&\max_{\Xi} \quad \sum_{t \in \mathcal{T}}\frac{1}{(1+r)^{t-1}}\bigg(\sum_{s \in \mathcal{S}} \sum_{b \in \mathcal{B}}\sum_{k \in \Omega_{t, s,b}^{D}}\Psi\pi_{t,s,b}d_{t, s,k,b} \notag\\
    &- \sum_{s \in \mathcal{S}} \sum_{ b\in \mathcal{B}}\sum_{k \in \Omega_{t, s,b}^{G}}\Psi\pi_{t,s,b}g_{t, s,k,b}+\Phi_t-\Psi\sum_{l \in \mathcal{L}}(u_{t,l} K_{l}^{f i x}\notag\\
    &+K_{l}^{v a r} \sum_{j \in \mathcal{J}} b_{t,l, j}^{F} \bar{F}_{l, j})\bigg)\label{eq:4a}
\end{align}
Subject to:
\begin{align}
\eqref{eq:upperb}-\eqref{eq:upperg}
\end{align}
\vspace{-18pt}
\begin{align}   
\eqref{eq:lowerb}-\eqref{eq:loweri}
\end{align}
\vspace{-18pt}
\begin{align}   
\eqref{eq:WSMdual}-\eqref{eq:strong_duality}
\label{eq:4d} \end{align}
\label{eq:4}
\end{subequations}
\end{fleqn}
where the variable array of the single-level problem~\eqref{eq:4} is described as $\Xi = \{u_{t,l},b_{t,l, j}^{F},\Phi_t,d_{t,s,k,b} ,g_{t,s,k,b},f_{t,s,l},\theta_{t,s,b},$ $[\pi _{t,s,b}]\}$. Notice that this single-level is a non-linear integer problem and next section will discuss the methods to remove non-linearities.

\subsection{Linearization}\label{section3c}
There exists two forms of non-linear terms in the single-level problem~\eqref{eq:4}:
\begin{enumerate}
\item the products $\pi_{t,s,b}d_{t,s, k, b}$ and $\pi_{t,s,b}g_{t,s, k, b}$ involving the WSM prices and clearing quantities in equations~\eqref{eq:upper_obj}, \eqref{eq:upperb} and~\eqref{eq:upperc}.
\item  the products $b^F_{\hat{t},l,j}\mu_{t,s, l}^{\max }$ and $b^F_{\hat{t},l,j}\mu_{t, s,l}^{\min }$ involving the binary variables $b^F_{\hat{t},l,j}$ and the continuous dual variables $\mu_{t,s, l}^{\max}$ and $\mu_{t,s, l}^{\min}$ in equation \eqref{eq:strong_duality}.
\end{enumerate}

To remove the first type of nonlinear terms, we exploit the definition of $\pi_{t,s,b}$ in equations~\eqref{eq:WSMduala}-\eqref{eq:WSMdualb}, we have 
\begin{flalign}
&\pi_{t, s,b} g_{t,s, k, b}= ( c_{t, s,k, b}^{g} +\varphi^{G,\max}_{t,s,k,b}-\varphi^{G,\min}_{t,s,k,b})g_{t, s,k, b} \label{eq:linearise1a} \\
&\pi_{t, s,b} d_{t, s,k, b} = ( c_{t, s,k, b}^{d} -\varphi^{D,\max}_{t,s,k,b}+\varphi^{D,\min}_{t,s,k,b})d_{t, s,k, b}  \label{eq:linearise1b}
\end{flalign}
Furthermore, the strong duality property enforced in equation~\eqref{eq:strong_duality} guarantees that all complementary slackness conditions hold. Therefore, we have 
\begin{equation}
\begin{aligned}
&\varphi^{G,\max}_{t,s,k,b}g_{t, s,k, b} = \varphi^{G,\max}_{t,s,k,b}g_{t, s,k, b}^{\max}\\ &\varphi^{G,\min}_{t,s,k,b}g_{t, s,k, b} = \varphi^{G,\min}_{t,s,k,b}g_{t, s,k, b}^{\min}\\
&\varphi^{D,\max}_{t,s,k,b}d_{t, s,k, b} = \varphi^{D,\max}_{t,s,k,b}d_{t, s,k, b}^{\max}\\
&\varphi^{D,\min}_{t,s,k,b}d_{t, s,k, b} = \varphi^{D,\min}_{t,s,k,b}d_{t, s,k, b}^{\min}
\end{aligned}  \label{eq:linearise1c}
\end{equation}
Then the terms $\pi_{t, s,b} g_{t,s, k, b}$ and $\pi_{t, s,b} d_{t,s, k, b}$ can be linearized as follows
\begin{flalign}
&\pi_{t, s,b} g_{t,s, k, b} \notag\\
&=    c_{t, s, k, b}^{g}g_{t, s, k, b} +\varphi^{G,\max}_{t,s,k,b}g_{t, s, k, b}^{\max}-\varphi^{G,\min}_{t,s,k,b}g_{t, s,k, b}^{\min} \label{eq:linearise1d}\\
&\pi_{t, s,b} d_{t, s,k, b}\notag\\
&=   c_{t,s, k, b}^{d}d_{t, s,k, b} -\varphi^{D,\max}_{t,s,k,b}d_{t, s,k, b}^{\max}+\varphi^{D,\min}_{t,s,k,b}d_{t,s, k, b}^{\min}\label{eq:linearise1e}
\end{flalign}
For the second type of nonlinearities, the big-M method is utilized \cite{SAVELLI2020113979}. The following constraints define two auxiliary variables $y^{\max}_{\hat{t},t,s,l,j}$ and $y^{\min}_{\hat{t},t,s,l,j}$ that are used to replace $b^F_{\hat{t},l,j}\mu_{t,s, l}^{\max }$ and $b^F_{\hat{t},l,j}\mu_{t, s,l}^{\min }$, respectively:
\begin{equation}
\begin{aligned}
& 0 \leq y^{\max}_{\hat{t},t,s,l,j}\leq M b_{\hat{t},l, j}^{F}\\
& 0 \leq \mu_{t,s,l}^{\max }-y^{\max}_{\hat{t},t,s,l,j} \leq M(1-b_{\hat{t},l, j}^{F})\end{aligned}  \label{eq:linearise3a}
\end{equation}

\begin{equation}
\begin{aligned}
& 0 \leq y^{\min}_{\hat{t},t,s,l,j} \leq M b_{\hat{t},l, j}^{F} \\
&0 \leq \mu_{t,s,l}^{\min }-y^{\min}_{\hat{t},t,s,l,j} \leq M(1-b_{\hat{t},l, j}^{F})  \end{aligned}  \label{eq:linearise3b}
\end{equation}
where constraints~\eqref{eq:linearise3a} and~\eqref{eq:linearise3b} are defined $\forall \hat{t} \leq t$, $\forall t \in \mathcal{T}$, $\forall s\in\mathcal{S}$, $\forall l\in\mathcal{L}$, $\forall j\in\mathcal{J}$. The parameter $M$ is the maximum allowed bid price in the European WSM, which equals \EUR3000/MWh.

\subsection{Final MILP problem}\label{section3d}
This section reports the full MILP model after all non-linearities have been removed as discussed in Section~\ref{section3c}.
\begin{fleqn}
\begin{subequations}
\begin{align}
&\max_{\Xi'} \sum_{t \in \mathcal{T}} \frac{1}{(1+r)^{t-1}}\bigg(\Psi\sum_{s \in \mathcal{S}}\sum_{b \in \mathcal{B}}\sum_{k \in \Omega_{t,s, b}^{D}}\big(c_{t, s,k, b}^{d}d_{t, s,k, b} \notag\\
&-\varphi^{D,\max}_{t,s,k,b}d_{t,s, k, b}^{\max} +\varphi^{D,\min}_{t,s,k,b}d_{t,s, k, b}^{\min} \big)\notag\\
&- \Psi\sum_{s \in \mathcal{S}}\sum_{ b\in \mathcal{B}}\sum_{k \in \Omega_{t, s,b}^{G}} \big( c_{t,s, k, b}^{g}g_{t,s, k, b}+\varphi^{G,\max}_{t,s,k,b}g_{t,s, k, b}^{\max}\notag\\
&-\varphi^{G,\min}_{t,s,k,b}g_{t,s, k, b}^{\min}\big)+\Phi_t\notag\\
&-\Psi\sum_{l \in \mathcal{L}}\big(u_{t,l} K_{l}^{f i x}+K_{l}^{v a r} \sum_{j \in \mathcal{J}} b_{t,l, j}^{F} \bar{F}_{l, j}\big)\bigg) \label{eq:finala}
\end{align}
Subject to:
\begin{align}
& S_{t}^{L} +S_{t}^{G}\notag\\
&=\sum_{s \in \mathcal{S}}\sum_{b \in \mathcal{B}}\bigg(\sum_{k \in \Omega_{t,s, b}^{D}}\big(\varphi^{D,\max}_{t,s,k,b}d_{t, s,k, b}^{\max}  -\varphi^{D,\min}_{t,s,k,b}d_{t, s,k, b}^{\min}\big)\notag\\
&+\sum_{k \in \Omega_{t, s,b}^{G}}\big(\varphi^{G,\max}_{t,k,b}g_{t, k, b}^{\max}-\varphi^{G,\min}_{t,k,b}g_{t, k, b}^{\min}\big)\bigg),\forall t\in\mathcal{T} \label{eq:finalb} \end{align}
\vspace{-15pt}
\begin{align}   
\eqref{eq:upperd}-\eqref{eq:upperg}
\label{eq:finalc} \end{align}
\vspace{-18pt}
\begin{align}   
\eqref{eq:lowerb}-\eqref{eq:loweri}
\label{eq:finald} \end{align}
\vspace{-18pt}
\begin{align}   
\eqref{eq:WSMdual}
 \label{eq:finale} \end{align}
\vspace{-18pt}
\begin{align}
&\sum_{t \in \mathcal{T}} \sum_{s \in \mathcal{S}}\sum_{b \in \mathcal{B}}\bigg(\sum_{k \in \Omega_{t, s,b}^{D}} c_{t,s, k,b}^{d} d_{t,s, k,b}-\sum_{k \in \Omega_{t, s,b}^{G}} c_{t,s,k,b}^{g} g_{t,s, k,b}\bigg)\notag\\
&=\sum_{t \in \mathcal{T}}\sum_{s \in \mathcal{S}} \bigg( \sum_{ b\in \mathcal{B}}\sum_{k \in \Omega_{t,s, b}^{G}}\big(\varphi^{G,\max}_{t,s,k,b}g^{\max}_{t,s,k,b} -\varphi^{G,\min}_{t,s,k,b}g^{\min}_{t,s,k,b}\big)\notag\\
&+\sum_{ b\in \mathcal{B}}\sum_{k \in \Omega_{t, s,b}^{D}}\big(\varphi^{D,\max}_{t,s,k,b}d^{\max}_{t,s,k,b}-\varphi^{D,\min}_{t,s,k,b}d^{\min}_{t,s,k,b}\big)\label{eq:finalf}\\ 
&+ \sum_{l \in \mathcal{L}}\mathcal{F}^{0}_{l}\big(\mu_{t,s, l}^{\max }+\mu_{t, s,l}^{\min }\big)+\sum_{b \in \mathcal{B}}\theta^{max}\big(\xi_{t,s, b}^{\max }+\xi_{t,s, b}^{\min }\big)\bigg)\notag\\
&+\sum_{t\in\{\mathcal{T}\backslash 1\}}\sum_{\hat{t}\leq t}  \sum_{s \in \mathcal{S}} \sum_{l \in \mathcal{L}}\sum_{j\in\mathcal{J}}\overline{F}_{l,j}\big(y^{\max}_{\hat{t},t,s,l,j} + y^{\min}_{\hat{t},t,s,l,j}\big)\notag
\end{align}
\vspace{-15pt}
\begin{align}   
\eqref{eq:linearise3a}-\eqref{eq:linearise3b}
\label{finalsocph} \end{align}
 \label{eq:finalsocp}
\end{subequations}
\end{fleqn}
where the variable array of the final MILP problem \eqref{eq:finalsocp} is $\Xi' = \{u_{t,l},b_{t,l, j}^{F},\Phi_t,d_{t,s,k,b} ,g_{t,s,k,b},f_{t,s,l},\theta_{t,s,b},y^{\max}_{\hat{t},t,s,l,j},$ $y^{\min}_{\hat{t},t,s,l,j},[\pi _{t,s,b}]\}$.
\section{Case Study}\label{section4}
Two case studies are presented to investigate the impact of the incentive tuning parameter on the Transco's profits, line expansion decisions, social welfare and the benefits for market participants. The first case study is a 2-node transmission network and the second one is the Garver's 6-node system. For simplicity and clarity of presentation, we assume that each investment planning period represents one year and includes one operation period $\mathcal{S}=\{1\}$. Therefore, the number of operation periods in one investment period $\Psi$ is $24\times365=8760$. These two case studies are implemented using PYOMO \cite{bynum2021pyomo} and CPLEX 22.1.0.0 and solved using a 13th Gen Intel(R) Core(TM) i7-13700K CPU @ 3.40GHz, 16 Core(s) with 32 Gb RAM. 
\subsection{2-node case study}
\begin{figure}[ht!]
    \centering    
    \includegraphics[width=\columnwidth]{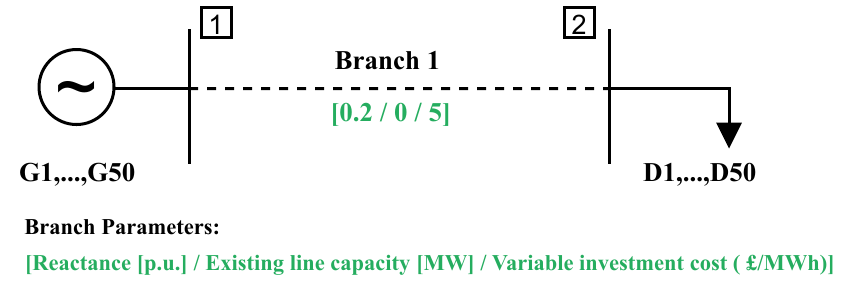}
    \caption{Topology of the 2-node transmission network.}
    \label{fig:2nodediagram}
\end{figure}

\begin{figure}[tp]
    \centering    
    \includegraphics[width=\columnwidth]{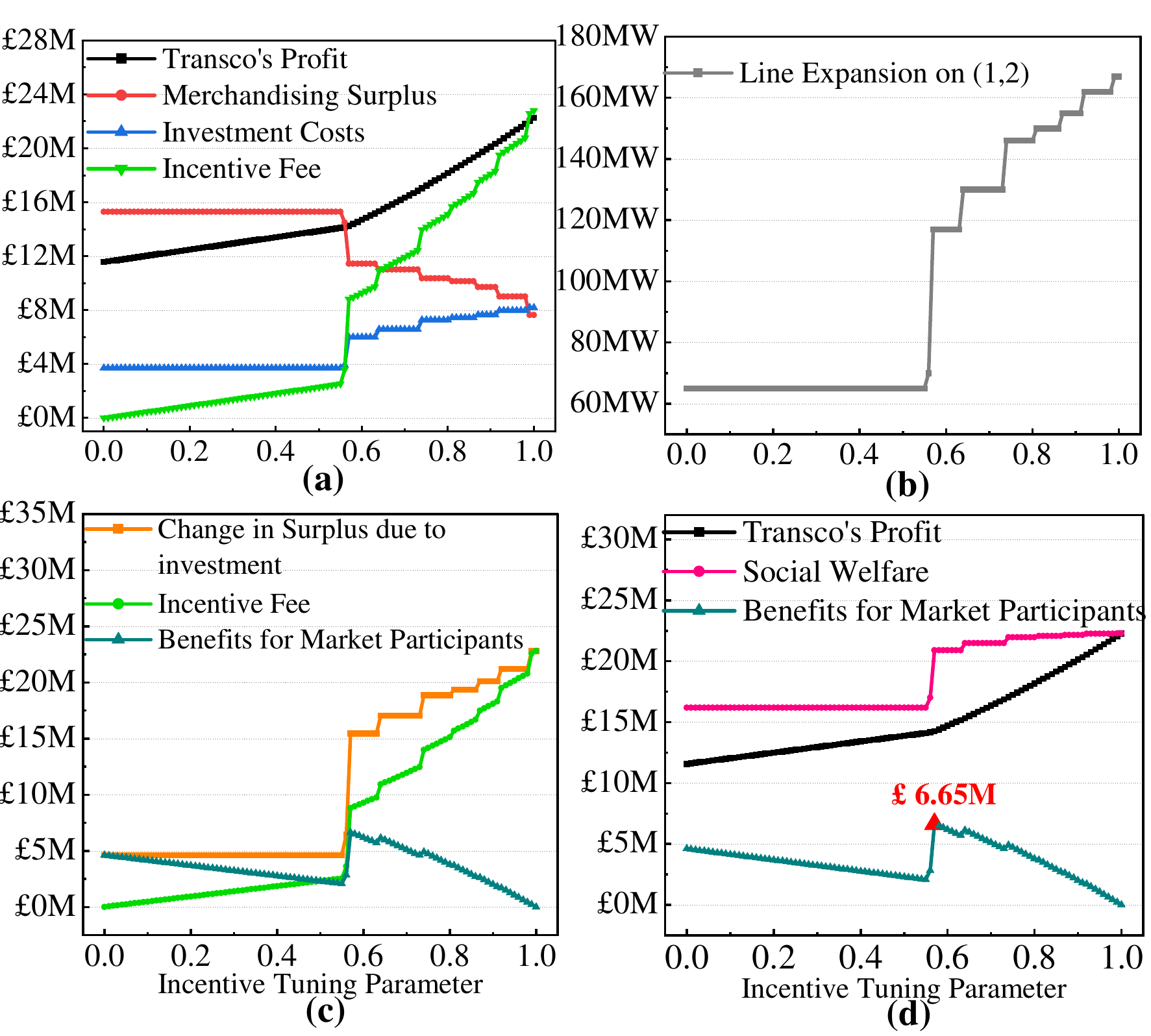}
    \caption{Results for the 2-node transmission network; (a) illustrates the Transco's profit, merchandising surplus, investment costs, and the incentive fee under different incentive tuning parameters. (b) shows the investment decision on line (1,2), while (c) presents the change in generator and consumer surplus due to investment (see the definition in Section~\ref{section2a}), the incentive fee, and the benefits for market participants. The incentive fee is computed by multiplying $\kappa$ with the change in surplus resulting from the investment (indicated by the orange line), with the remaining portion representing the benefits for market participants. (d) provides a comparison among the private interests (the Transco's profit), public interests (social welfare), and market participants' interests (benefits for market participants).\color{black}}
    \label{fig:2node1}
\end{figure}

The first case study is based on a 2-node transmission network in which node 2 has 50 consumers and node 1 has 50 generators, as shown in Fig.~\ref{fig:2nodediagram}. Bid prices for consumers $c^{d}_{t,s,k,b}$ and generators $c^{d}_{t,s,k,b}$ are randomly generated from normal distributions with mean prices of $\pounds$50/MWh and $\pounds$40/MWh, respectively, with the same standard deviations of $\pounds$10/MWh. The limits of generations at node 1 and demands at node 2 are generated from a uniform distribution ranging from zero to 10 MW. The line reactance is set to 0.2 p.u. and the existing line capacity between two nodes is zero. The variable investment cost $K_l^{var}$ is $\pounds$5/MWh and the fixed cost $K_l^{fix}$ is $\pounds$100/h. The set of lumpy capacity expansions is defined as $\mathcal{F}_l=\{1,2,...,400\}$MW. The operational timescale of the planning problem includes two investment years, i.e., $\mathcal{T} = \{1,2\}$. Different values of the incentive tuning parameter $\kappa = \{0,0.01,0.02,...,1\}$ were considered over two years, and the results are shown in Fig.~\ref{fig:2node1}. Table~\ref{table1} reports the highlighted results of the 2-node transmission network. The average simulation time for the problems was 1.3 seconds.

\color{black}

\begin{table}[ht]
\centering
\begin{threeparttable}
\caption{Statistics on Private, Public, and Market Participant Interests of the Two-Node Transmission Network Case Study}
\label{table1}
\begin{tabular}{@{}l>{\centering\arraybackslash}p{2cm}>{\centering\arraybackslash}p{2cm}>{\centering\arraybackslash}p{2.5cm}@{}}
\hline
 & \textbf{The Transco's profit [M$\pounds$]} & \textbf{Social Welfare [M$\pounds$]} & \textbf{Benefits for Market Participants [M$\pounds$]} \\
\hline
$\kappa = 1$\tnote{*} & 22.26 & 22.26 & 0 \\
$\kappa = 0.57$\tnote{**} & 14.25 & 20.89 & 6.65 \\
$\kappa = 0$\tnote{***} & 11.58 & 16.18 & 4.60 \\
\hline
\end{tabular}
\begin{tablenotes}
\item[*] $\kappa = 1$: This setting corresponds to the original H-R-G-V mechanism, under which the entire surplus increase is awarded to the Transco, leaving market participants without any share of the increased surplus.
\item[**] $\kappa = 0.57$: This value represents a modification to the H-R-G-V mechanism, designed to optimize the benefits received by market participants, thereby indicating a balanced distribution of surplus increases that favors market participants.
\item[***] $\kappa = 0$: This scenario, another modification of the H-R-G-V mechanism, ensures that market participants receive the entire surplus increase, with the Transco not receiving any portion of the surplus enhancement.
\end{tablenotes}
\end{threeparttable}
\end{table}

\color{black}
The decision to expand the line is influenced by the value of the incentive tuning parameter, $\kappa$, as depicted in Fig.~\ref{fig:2node1}(b). Notably, the line expansion decision reaches its minimum at 65 MW when $\kappa$ falls within the range of $\{0,...,0.55\}$. This can be attributed to the Transco not receiving sufficient incentives to perform investments, and the benefit it can obtain from a higher merchandising surplus resulting from line congestion. Consequently, the social welfare (the pink line in Fig.~\ref{fig:2node1}(d)) and the change in generator and consumer surplus due to investment (the orange line in Fig.~\ref{fig:2node1}(c), see the definition in Section~\ref{section2a}) \color{black}reach their lowest values at $\pounds$16.18M and $\pounds$4.60M, respectively. Within this region, there is a consistent increase in the incentive fee and a corresponding decrease in the benefits for market participants as the incentive parameter increases. This trend arises from the fact that the change in surplus due to investment remains the same, and the portion of the change in surplus multiplied by $\kappa$ belongs to the Transco, while the remaining benefits are received by market participants. 

The line expansion decision rises to 70 MW when $\kappa$ equals 0.56 and undergoes a significant increase, reaching 117 MW within the range of $\kappa\in\{0.57,...,0.63\}$. As illustrated in Fig.~\ref{fig:2node1}(d), the benefits for market participants reach their peak value of $\pounds$6.65M at $\kappa=0.57$, signifying that this value represents an \emph{optimal} decision from the market participants' perspective. As shown in Table~\ref{table1}, increasing $\kappa$ from 0 to 0.57 results in a 29$\%$ increase in social welfare, rising from $\pounds$16.18M to $\pounds$20.89M, while the Transco's profit also experiences a significant increase of nearly 33.79$\%$, ascending from $\pounds$11.58M to $\pounds$14.25M.

With a further increase in the tuning parameter, the amount of line expansion decision and investment costs exhibits a fluctuating pattern, characterized by intermittent intervals of constancy but an overall upward trend. This behavior aligns with the observed trends in social welfare and the changes in generator and consumer surplus resulting from investment. Conversely, the benefits for market participants exhibit a noncontinuous downward trajectory, declining from $\pounds$6.65M to 0 as the incentive tuning parameter increases from 0.57 to 1. This divergence arises due to a greater proportion of benefits being allocated to the Transco rather than market participants for higher values of $\kappa$.

When $\kappa = 1$, the line (1,2) is built with a capacity of $\bar{F}_1 = 167$ MW, resulting in the maximum welfare of $\pounds$22.26M. The total investment cost in this scenario amounts to $\pounds$8.19M, with the Transco earning the highest profit of $\pounds$22.26M, including an incentive fee of $\pounds$22.80M. It is evident that at $\kappa = 1$, the proposed scheme replicates the standard H-R-G-V scheme, where the Transco captures the entire surplus change. 

In Fig.~\ref{fig:2node1}(d) and Table~\ref{table1}, when the value of $\kappa$ is reduced from 1 to 0.57, there is only a slight decrease of 6.15$\%$ in social welfare, from $\pounds$22.26M to $\pounds$20.89M. However, the benefits for market participants increase substantially, rising from 0 to $\pounds$6.65M. Conversely, the Transco's profit drops by nearly 36$\%$, from $\pounds$22.26M to $\pounds$14.25M. Hence, by decreasing the value of $\kappa$, the Transco's profits can be constrained, while still achieving a modest level of social welfare.

\subsection{6-node case study}
\begin{figure}[tp]
    \centering    
    \includegraphics[width=\columnwidth]{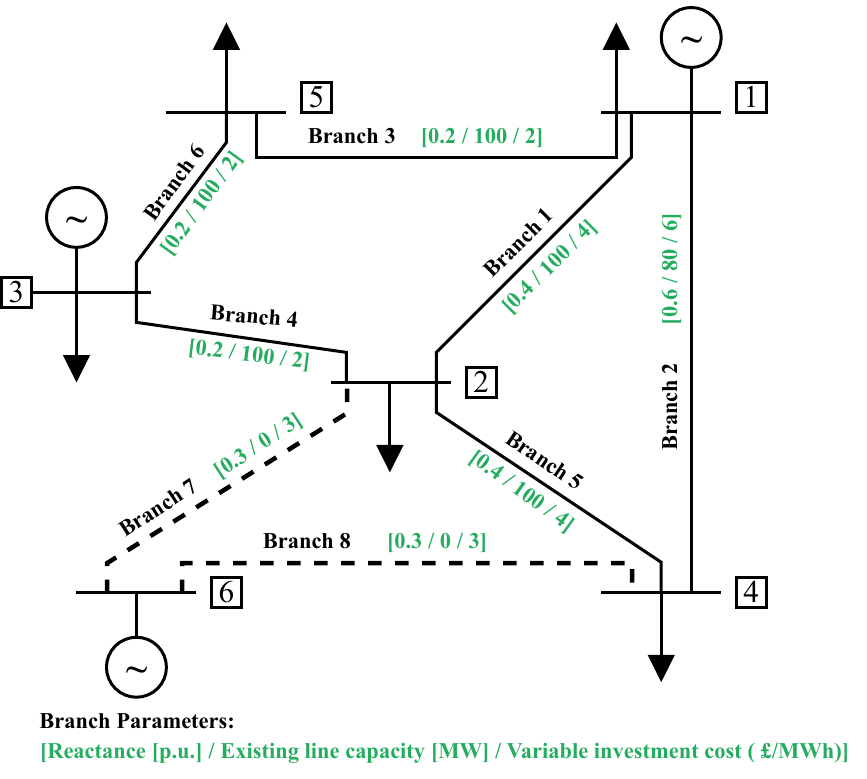}
    \caption{Topology of the Garver's 6-node transmission network.}
    \label{fig:6nodediagram}
\end{figure}

\begin{figure}[tp]
    \centering    \includegraphics[width=\columnwidth]{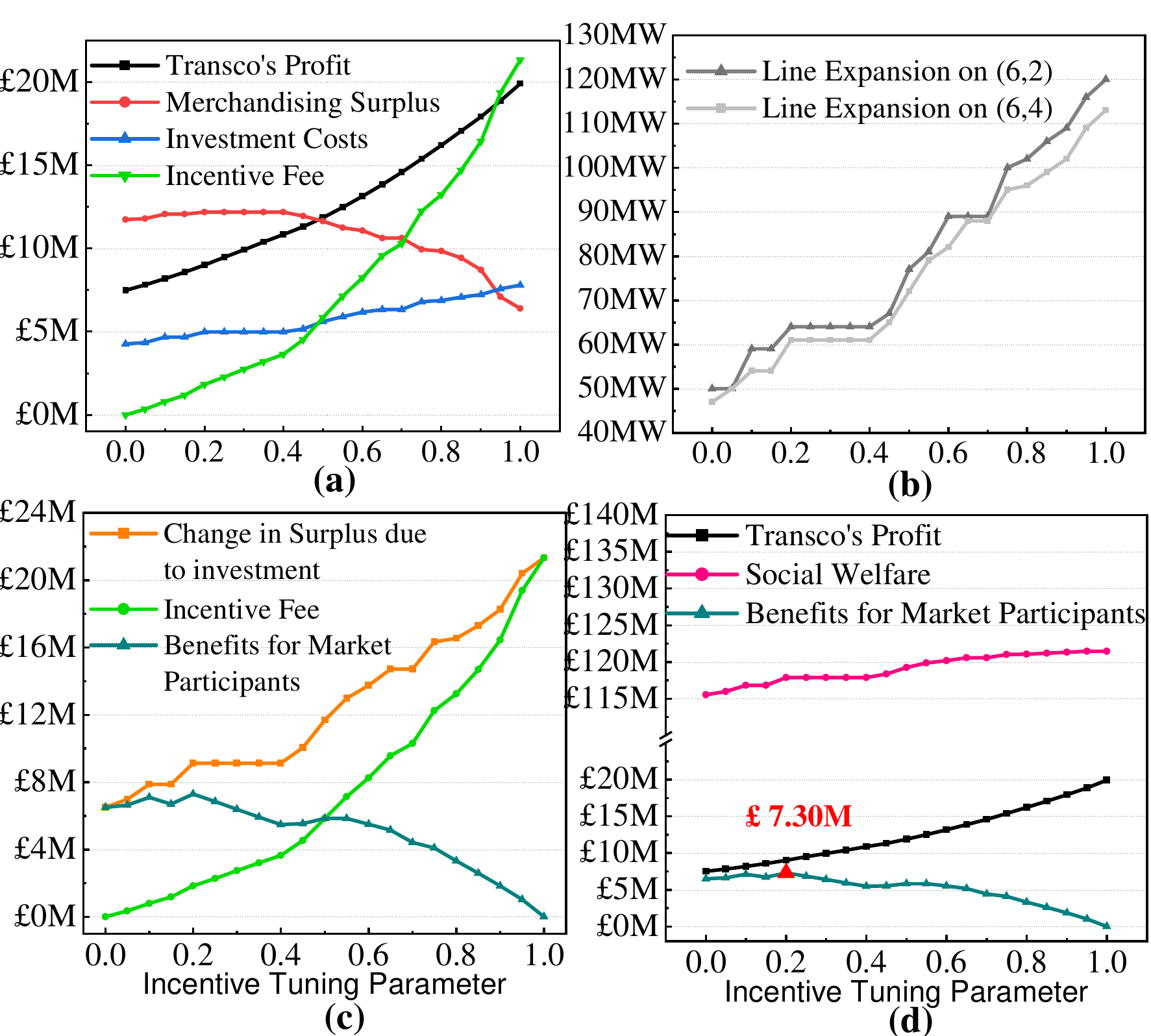}
    \caption{Results for the Garver's 6-node transmission network.}
    \label{fig:6node1}
\end{figure}

The case study investigates the investment planning problem in a more complex setting using the modified Garver's 6-node transmission network \cite{SAVELLI2020113979,Pavel,Mazaheri}. This network consists of 6 nodes and 8 lines, denoted as $\mathcal{B} = \{1,...,6\}$ and $\mathcal{L}=\{1,...,8\}$, as shown in Fig.~\ref{fig:6nodediagram}. The investment planning problem focuses on the expansion of existing lines (branches 1 to 6) and the construction of new lines (branches 7 and 8). In this study, we adopt a case study setting similar to that in \cite{SAVELLI2020113979}, with nodes $\{1,3,6\}$ each having 1,000 generators and nodes $\{1,2,3,4,5\}$ each accommodating 1,000 consumers to accurately depict elastic market curves. This design of generators and consumers at specific nodes follows the methodology used in research by \cite{Pavel,Mazaheri}. Bid prices for consumers $c^{d}_{t,s,k,b}$ and generators $c^{d}_{t,s,k,b}$ are randomly generated from normal distributions with mean prices of $\pounds$50/MWh with the standard deviations of $\pounds$10/MWh. The limits of demands at nodes $\{1,2,3,4,5\}$ and generators at nodes 1 and 3 are generated from a uniform distribution, ranging from zero to 500 kW \cite{SAVELLI2020113979}. For node 6, we establish the supply limit with random numbers between zero and 1 MW, also using a uniform distribution, to create an interesting case study. The annual growth rate of load is assumed to be 5$\%$. The operational timescale of the planning problem includes two investment periods, i.e., $\mathcal{T} = \{1,2\}$ and the discount rate $r$ is set to 1$\%$. The findings of the 6-node system are presented in Fig.~\ref{fig:6node1} for different values of $\kappa = \{0,0.05,0.1,...,1\}$. Table~\ref{table2} reports the highlighted results of the Garver's 6-node transmission network. The average simulation time for the problems was 2,515 seconds.

\begin{table}[ht]
\centering
\begin{threeparttable}
\caption{Statistics on Private, Public, and Market Participant Interests of the Garver's Six-Node Transmission Network Case Study}
\label{table2}
\begin{tabular}{@{}l>{\centering\arraybackslash}p{2cm}>{\centering\arraybackslash}p{2cm}>{\centering\arraybackslash}p{2.5cm}@{}}
\hline
 & \textbf{The Transco's profit [M$\pounds$]} & \textbf{Social Welfare [M$\pounds$]} & \textbf{Benefits for Market Participants [M$\pounds$]} \\
\hline
$\kappa = 1$ & 19.91 & 121.45 & 0 \\
$\kappa = 0.20$  & 9.01 & 117.85 & 7.30 \\
$\kappa = 0$ & 7.48 & 115.51 & 6.49\\
\hline
\end{tabular}
\end{threeparttable}
\end{table}

The line expansion decisions on line (6,2) and (6,4) are depicted in Fig.~\ref{fig:6node1}(b). There is a general upward trend in the line expansion decision with the increasing tuning parameters from 0 to 1, where the expansion amount increase from 50 MW and 47 MW to 120 MW and 113 MW for line (6,2) and line (6,4), respectively. Similar to the results for the 2-node transmission network, we also notice the increase in investment costs and the change in surplus and social welfare aligns with the line expansion results pattern, and these curves remains constant during some intervals. In addition, the Transco's profit displays a consistent upward trend with the increasing tuning parameter $\kappa$. This phenomenon can be attributed to two factors. Firstly, as the line expansion decision increases, the total surplus change resulting from the investment also experiences a significant increase. Secondly, the incentive fee received by the Transco demonstrates an increase, indicating that a larger proportion of the surplus change is allocated to the Transco as the tuning parameter increases.

Fig.~\ref{fig:6node1}(d) provides a reference for regulators in determining the optimal incentive tuning parameter. As expected, the social welfare exhibits a general increasing trend as the tuning parameter increases. Notably, when $\kappa=0.20$, market participants receive the maximum benefit, amounting to $\pounds$7.30M while the social welfare only decreases by $3\%$ from $\pounds$121.45M (when $\kappa=1$) to $\pounds$117.85M (when $\kappa=0.20$) . This suggests that from the perspective of market participants, $\kappa = 0.2$ represents the \emph{optimal} decision. Similar to the 2-node case study, we observe that when $\kappa = 1$, the incentive fee calculated based on the total change in surplus due to investment emerges as the optimal strategy for optimizing social welfare and the Transco's profit. This can be attributed to the fact that both the Transco's profit and social welfare reach their peak levels under this scenario.

\color{black}
\section{Conclusion}\label{section5}

This paper presents an extension to the H-R-G-V mechanism by introducing the incentive tuning parameter to appropriately incentivize the Transco to invest in the transmission network considering the increasing price-elasticity of demand. In the proposed mechanism, the regulated incentive fee is directly linked to the total economic benefits generated by the investment and the choice of this tuning parameter. This approach benefits market participants by allowing them to capture the economic benefits resulting from the investment.

The application of the mechanism is tested through two case studies, and results show that the Transco is performing social welfare maximizing investment when the Transco receives all generator and consumer surplus increase due to investment (i.e., $\kappa=1$ for both cases). Consumers and generators receive the maximum economic benefit with reduced tuning parameter (i.e., $\kappa=0.57$ for the 2-node case study and $\kappa=0.20$ for the 6-node case study). The allocation of these benefits to market participants contributes to a more fair and equitable distribution of the investment's positive outcomes. Moreover, the proposed regulatory mechanism can significantly increase the benefits for consumers and generators from transmission network investments while the impact on overall social welfare remains relatively modest.

Future work could incorporate the uncertainties associated with load profile, environmental policies and renewable generation, focusing on developing a stochastic model considering multiple scenarios. Another interesting research area is to integrate the reliability requirements for the secure operation of the network.

\color{black}
\appendix
This appendix presents dual constraints and the strong duality condition of the lower-level problem.
\subsection{Lower-level dual constraints}\label{section3a}
The dual constraints of the lower-level WSM clearing problem are presented in~\eqref{eq:WSMdual}.
\begin{fleqn}
\begin{subequations}
\begin{align}
& -\pi_{t,s,b}+\varphi^{G,\max}_{t,s,k,b}-\varphi^{G,\min}_{t,s,k,b}= - c^{g}_{t,s,k,b}, \notag\\
& \forall  t\in\mathcal{T}, \forall s\in \mathcal{S},\forall k \in \Omega_{t, s,b}^{G},b\in\mathcal{B} \hspace{25pt}[g^p_{t,s,k,b}\in\mathbb{R}] \label{eq:WSMduala} \end{align}
\vspace{-10pt}
\begin{align}
&\pi_{t,s,b}+\varphi^{D,\max}_{t,s,k,b} -\varphi^{D,\min}_{t,s,k,b}=  c^{d}_{t,s,k,b},\notag\\
& \forall t\in\mathcal{T}, \forall s\in \mathcal{S},\forall k \in \Omega_{t, s,b}^{D},b\in\mathcal{B} \hspace{25pt}[d^{p}_{t,s,k,b}\in\mathbb{R}] \label{eq:WSMdualb} \end{align}
\vspace{-10pt}
\begin{align}
&\sum_{b\in\mathcal{B}}S_{l,b}\pi_{t,s,b}-\sum_{b\in\mathcal{B}}R_{l,b}\pi_{t,s,b}+\gamma_{t,s,l} + \mu^{max}_{t,s,l}-\mu^{min}_{t,s,l}=0, \notag\\
& \forall t\in\mathcal{T},\forall s\in \mathcal{S}, \forall l \in \mathcal{L} \hspace{80pt}[f_{t,s,l}\in\mathbb{R}]\label{eq:WSMdualc} \end{align}
\vspace{-10pt}
\begin{align}
&-B_{l}\sum_{l\in\mathcal{L}}S_{l,b}\gamma_{t,s,m}+B_{l}\sum_{l\in\mathcal{L}}R_{l,b}\gamma_{t,s,l} + \xi^{max}_{t,s,b}-\xi^{min}_{t,s,b}=0, \notag\\
& \forall  t\in\mathcal{T}, \forall s\in \mathcal{S},\forall b \neq 1  \hspace{80pt}[\theta_{t,s,b}\in\mathbb{R}] \label{eq:WSMduald} \end{align}
\vspace{-10pt}
\begin{align}
&-B_{l}\sum_{l\in\mathcal{L}}S_{l,b}\gamma_{t,s,m}+B_{l}\sum_{l\in\mathcal{L}}R_{l,b}\gamma_{t,s,l} + \chi_{t,s,b}=0, \notag \\
&\forall  t\in\mathcal{T},\forall s\in \mathcal{S}, b = 1 \hspace{85pt}[\theta_{t,s,1}\in\mathbb{R}] \label{eq:WSMduale} \end{align}
\label{eq:WSMdual}
\end{subequations}
\end{fleqn}
\subsection{Lower-level Strong Duality Condition}\label{section3b}
The strong duality condition~\eqref{eq:strong_duality} requires the equivalence between the primal and dual objective values, as shown below:
\begin{align}
            &\sum_{t \in \mathcal{T}} \sum_{s \in \mathcal{S}}\sum_{b \in \mathcal{B}}\bigg(\sum_{k \in \Omega_{t, s,b}^{D}} c_{t,s, k,b}^{d} d_{t,s, k,b}-\sum_{k \in \Omega_{t, s,b}^{G}} c_{t,s,k,b}^{g} g_{t,s, k,b}\bigg)\notag\\
            &=\sum_{t \in \mathcal{T}}\sum_{s \in \mathcal{S}} \bigg( \sum_{ b\in \mathcal{B}}\sum_{k \in \Omega_{t,s, b}^{G}}\big(\varphi^{G,\max}_{t,s,k,b}g^{\max}_{t,s,k,b} -\varphi^{G,\min}_{t,s,k,b}g^{\min}_{t,s,k,b}\big)\notag\\
            &+\sum_{ b\in \mathcal{B}}\sum_{k \in \Omega_{t, s,b}^{D}}\big(\varphi^{D,\max}_{t,s,k,b}d^{\max}_{t,s,k,b}-\varphi^{D,\min}_{t,s,k,b}d^{\min}_{t,s,k,b}\big)\label{eq:strong_duality}\\ 
            &+ \sum_{l \in \mathcal{L}}\mathcal{F}^{0}_{l}\big(\mu_{t,s, l}^{\max }+\mu_{t, s,l}^{\min }\big)+\sum_{b \in \mathcal{B}}\theta^{\max}_b\big(\xi_{t,s, b}^{\max }+\xi_{t,s, b}^{\min }\big)\bigg)\notag\\
            &+\sum_{t\in\{\mathcal{T}\backslash 1\}}\sum_{\hat{t}\leq t}  \sum_{s \in \mathcal{S}} \sum_{l \in \mathcal{L}}\sum_{j\in\mathcal{J}}\overline{F}_{l,j}\big(b^F_{\hat{t},l,j}\mu_{t,s, l}^{\max }+b^F_{\hat{t},l,j}\mu_{t, s,l}^{\min }\big)\notag
            \end{align} 

\nomenclature[S]{$\mathcal{L}$}{Set of transmission lines}
\nomenclature[S]{$\mathcal{J}$}{Set of lumpy capacity indices}
\nomenclature[S]{$\mathcal{B}$}{Set of transmission network nodes}
\nomenclature[S]{$\mathcal{T}$}{Set of investment periods}
\nomenclature[S]{$\mathcal{S}$}{Set of operation periods}
\nomenclature[S]{$\Omega_{t,s,b}^{G}$}{ Set of generators at investment period $t$, operation period $s$ in node $b$}
\nomenclature[S]{$\Omega_{t,s,b}^{D}$}{ Set of consumers at investment period $t$, operation period $s$ in node $b$}
\nomenclature[S]{$\overline{\mathcal{F}}_{l,j}$}{Lumpy capacity expansion for line $l$, with $\overline{\mathcal{F}}_l=\bigcup_{j \in \mathcal{J}} \bar{F}_{l, j}$}
\nomenclature[P]{$r$}{The discount rate}
\nomenclature[P]{$\kappa$}{Incentive tuning parameter, $\kappa \in [0\%, 100\%]$}
\nomenclature[P]{$\Psi$}{Number of operation periods in one year}
\nomenclature[P]{$\theta^{max}_b$}{Maximum voltage angle at node $b$ (rad)}
\nomenclature[P]{$c^d_{t,s,k,b}$}{Demand bid price (willingness-to-pay) for consumers $k$ in node $b$ at year $t$ and period $s$ ($\pounds$/MWh)}
\nomenclature[P]{$c^p_{t,s,k,b}$}{Supply bid price (marginal cost) for generators $k$ in node $b$ at year $t$ and period $s$ ($\pounds$/MWh)}
\nomenclature[P]{$K_l^{fix}$}{Fixed cost of building or expanding line $l$ ($\pounds$/h)}
\nomenclature[P]{$K_l^{var}$}{Variable cost of building or expanding line $l$ ($\pounds$/MWh)}
\nomenclature[P]{$B_{l}$}{Susceptance of the transmission line $l$, $l \in \mathcal{L}$ (S)}
\nomenclature[P]{$S_{l,b}$}{Incidence matrix element of sending node $b$, line $l$, $b \in \mathcal{B}$, $l \in \mathcal{L}$}
\nomenclature[P]{$R_{l,b}$}{Incidence matrix element of receiving node $b$, line $l$, $b \in \mathcal{B}$, $l \in \mathcal{L}$}
\nomenclature[P]{$\mathcal{F}^{0}_{l}$}{ Existing capacity on the transmission line $l$, $l \in \mathcal{L}$ (MW)}
\nomenclature[P]{$d_{t,s, k,b}^{\min},d_{t,s, k,b}^{\max}$}{Minimum/maximum quantity of active power demanded by consumers $k$ at the investment period $t$, operation period $s$ in node $b$ (MW)}
\nomenclature[P]{$g_{t,s, k,b}^{\min},g_{t,s, k,b}^{\max}$}{Minimum/maximum quantity of active power produced by generators $k$ at the investment period $t$, operation period $s$ in node $b$ (MW)}
\nomenclature[V]{$u_{t,l}$}{Binary variable equal to one if line $l$ is expanded at the investment period $t$, and zero otherwise, $l \in \mathcal{L}_b$}
\nomenclature[V]{$b_{t, l, j}^{F}$}{ Binary variable equal to one if the lumpy investment in additional capacity $\bar{F}_{l, j}$ for line $l$ is made at the investment period $t$, and zero otherwise, $l \in \mathcal{L}_b$}
\nomenclature[V]{$\Phi_t$}{The incentive fee at the investment period $t$ (\pounds)}
\nomenclature[V]{$g_{t,s,k,b}$}{Allocated active power for generators $k$ at the investment period $t$, operation period $s$ in node $b$ (MW)}
\nomenclature[V]{$d_{t,s,k,b}$}{Allocated active power for consumers $k$ at the investment period $t$, operation period $s$ in node $b$ (MW)}
\nomenclature[V]{$f_{t,s,l}$}{Flow in the line $l$ at the investment period $t$, operation period $s$ (MW)}
\nomenclature[V]{$\theta_{t,s,b}$}{Voltage phase angle of transmission network node $b$ at the investment period $t$, operation period $s$ (rad)}
\nomenclature[V]{$\pi_{t,s,b}$}{WSM prices at transmission node $b$ at the investment period $t$, operation period $s$ (\pounds$/$MWh)}
\nomenclature[V]{$y^{\max}_{\hat{t},t,s,l,j}$}{Replace the product $b^F_{\hat{t},l,j}\mu_{t,s, l}^{\max}$}
\nomenclature[V]{$y^{\min}_{\hat{t},t,s,l,j}$}{Replace the product $b^F_{\hat{t},l,j}\mu_{t, s,l}^{\min}$}
\printnomenclature

\section*{Acknowledgement}
The authors would like to express our gratitude to D. Khastieva, M.R. Hesamzadeh, I. Vogelsang, and their colleagues for their influential work titled "Transmission Network Investment Using Incentive Regulation: A Disjunctive Programming Approach" \cite{Disjunctive}.\color{black}
\bibliographystyle{IEEEtran}
\bibliography{reference.bib}

\end{document}